\documentstyle[12pt]{article}

\font\tenmsb=msbm10 scaled\magstep1
\font\sevenmsb=msbm7 scaled\magstep1
\font\fivemsb=msbm5 scaled\magstep1
\newfam\msbfam
\textfont\msbfam=\tenmsb
\scriptfont\msbfam=\sevenmsb
\scriptscriptfont\msbfam=\fivemsb
\def\Bbb#1{{\fam\msbfam\relax#1}}

\newcommand{\cc}{{\Bbb C}}       % Field of complex numbers
\newcommand{\rr}{{\Bbb R}}       % Field of real numbers
\newcommand{\nn}{{\Bbb N}}       % Natural numbers
\newcommand{\T}{{\Bbb T}}        % Unit torus

\def\E{{\cal E}}

\def\supp{\mbox{\rm supp}}

\def\U{\Upsilon}
\def\L{\Lambda}
\def\g{\gamma}

\newtheorem{theorem}{Theorem}
\newtheorem{lemma}[theorem]{Lemma}
\newtheorem{proposition}[theorem]{Proposition}
\newtheorem{corolari}[theorem]{Corollary}
\newtheorem{definition}{Definition}
\newcommand{\demo}[1]{\par\smallskip \noindent {\bf {#1}:}}
\newcommand{\qed}{\ \hfill\mbox{$\clubsuit$}}

\title{On the existence of doubling measures with certain regularity
properties}

\author{Per Bylund\and Jaume
  Gudayol\thanks{Partially supported by MEC grant PB95-0956-C02-01
and CIRIT grant GRQ94-2014.}}

\date{xx March 1998}

\begin{document}

\maketitle

\begin{abstract}
Given a compact pseudo-metric space, we associate to it upper and lower
dimensions, depending only on the metric. Then we
construct a doubling measure for which the measure of a dilated ball is closely
related to these dimensions.
\end{abstract}

\section{Introduction}
Let $(X,\rho)$ be a compact complete metric space. Suppose that
$(X,\rho)$ is homogeneous. This means that there exists a so-called
doubling measure $\mu$ supported by $X$, i.e. there is a constant $c$ such
that for any $x\in X$ and any $R>0$
\begin{equation}
\mu(B(x,2R))\le c \mu(B(x,R)).\label{cal_posarho}
\end{equation}
In 1984, Dynkin (\cite{Dyn}) proved that for certain subsets $E$ of
the unit sphere $\T\subset \cc$ there exists a doubling measure on $E$.
In the same paper Dynkin conjectured that any compact $E\subset\rr^n$ is
homogeneous.
This conjecture was proved by Volberg and Konyagin (\cite{V-K}) by
using a dimension first defined in
\cite{L} under the name of uniform metric dimension,
in this paper denoted by $\Upsilon(E)$
(note that $\Upsilon(\rr^n)=n$ in the Euclidean case).
More precisely, Volberg and Konyagin proved that $(X,\rho)$ is homogeneous
if and only if there is some $\gamma<\infty$ such that any ball $B(x,kR)$
contains
at most $Ck^{\gamma}$ points separated from each other by a distance of at
least $R$.
The uniform metric dimension $\Upsilon(X,\rho)$ is then defined as the
infimum of such
$\gamma$.
Furthermore, given $\gamma<\infty$ in the condition above Volberg and Konyagin
proved that for any $s>\gamma$ there exists a measure $\mu$ such that,
for $0<R\le kR\le1$,
\begin{equation}
\mu(B(x,kR))\le C_1 k^s \mu(B(x,R))      \label{aixo_tambe}
\end{equation}
Clearly, any measure satisfying \ref{aixo_tambe} is a doubling
measure, and conversely, iterating \ref{cal_posarho} one gets
\ref{aixo_tambe} with $s=\log_2c$.
In particular, Volberg and Konyagin proved Dynkin's conjecture by showing
that on any
compact $E\subset
\rr^n$ there exists a measure $\mu$ satisfying \ref{aixo_tambe}
 with $s=n$ (in the
maximum metric).

In this paper we generalize the proof of Volberg and Konyagin to the
pseudo-metric case
by showing the existence of a measure $\mu$ not only
satisfying the upper bound condition \ref{aixo_tambe},
 but also the following analogous lower bound
condition: Suppose there is a $\delta\ge0$ such that any ball $B(x,kR)$
contains at
least $Ck^{\delta}$ points separated from each other by a distance of at
least $R$.
Then for any $t<\delta$ there exists a measure $\mu$ such that, for $0<R\le
kR\le1$,
\begin{equation}
C_2 k^t \mu(B(x,R))\le\mu(B(x,kR))      \label{otroquetal}
\end{equation}
Note that \ref{otroquetal} is trivially true for $t=0$.

Jonsson and Wallin in \cite{J-W} gave a thorough study of function spaces on
$s$-sets. By definition, an $s$-measure fulfils both \ref{aixo_tambe} and
\ref{otroquetal}
in the special case when $s=t$. An $s$-set is a set on which
there exists an $s$-measure, which then may be taken as
the $s$-dimensional Hausdorff measure. These sets are also called
Ahlfors-regular sets.

Measures satisfying both \ref{aixo_tambe}
and \ref{otroquetal} in the general case when $t\le s$ were
first considered
by  Jonsson (\cite{Jonsson})
when studying interpolation sets for Besov spaces on $\rr^n$.

The authors of this paper, independently of each other, also studied such
measures in
\cite{Per} and \cite{jo}. Each of these works contains the main result of this
paper,
in \cite{Per} formulated for Euclidean spaces and in \cite{jo} for metric
spaces.
However, in this paper our result has been restated in terms of
pseudo-metric spaces as we believe there are some applications to this more
general case.

\section{Definitions and statements of results}
In what follows we denote by $(X,d)$ a complete locally compact
pseudo-metric space.
We say that $d:X\times X \mapsto [0,+\infty)$ is a pseudo-metric on $X$ if the
following properties are fulfilled: \begin{enumerate}
\item $d(x,y)=0\iff x=y$
\item $d(x,y)=d(y,x),\quad\forall x,y\in X$
\item there is a constant $C_d$ such that $\forall x,y,z\in X$,
$$d(x,z)\le C_d(d(x,y)+d(y,z)).$$
\end{enumerate}
Note that $C_d\ge1$.

Given any ball $B(x,kR)$, $x\in X$ and $0<R\le kR$,
denote by $N(x,R,k)$ the maximum number of points in $B(x,kR)$
separated by a distance greater than or equal to $R$ from each other.

\begin{definition}
We will say that $(X,d)\in\Upsilon_\g$ if there exists
$C(\g)=C(X,d,\g)$ such that, for $kR\le1$,
$$
N(x,R,k)\le C(\g)k^\g.\eqno{(\Upsilon_\g)}
$$
Then we define the upper dimension $\Upsilon(X)$ as
$$
\Upsilon(X)=\inf\{\g,\,(X,d)\in\Upsilon_\g\}.
$$
\end{definition}
This dimension was first defined by Larman (\cite{L}) under the name of
uniform metric
dimension.

\begin{definition}
We will say that a probability measure $\mu$ lies in $U_\g=U_\g(X,d)$
if there exists $C(\g)$ such that, for $x\in X$ and $0<R\le kR\le1$,
$$
\quad\mu(B(x,kR))\le Ck^\g \mu(B(x,R))\eqno{(U_\g)}
$$
Then we define the  dimension $U(X)$ as
$$
U(X)=\inf\{\g,\,U_\g(X,d)\ne\emptyset\}.
$$
\end{definition}
Note that by taking $k=1/R$ in $(U_\g)$ one gets the weaker condition
$$
\mu(B(x,R))\ge C R^\g,\qquad x\in X,\quad 0<R\le1.\eqno{(U'_\g)}
$$
Also note that if $\mu\in U_\g$, for some $\g$, then $\supp\mu =X$.
As mentioned in the introduction, $\mu$ is doubling precisely
when $\mu\in U_\g$ for some $\g<\infty$. We will write ${\cal U}=\cup_\g
U_\g$ for
the set of all doubling measures on $X$.

Volberg and Konyagin in \cite{V-K} proved that $\Upsilon(X)\le U(X)$.
Furthermore they proved:
\begin{theorem}[Volberg-Konyagin] \label{VK4}
Let $(X,d)$ be a compact complete metric space. If $(X,d)\in\Upsilon_\g$ then
for every $\g'>\g$ there exists a measure $\mu\in U_{\g'}$. Consequently,
$\U(X,d)=U(X,d)$.
\end{theorem}

We will prove, as a part of Theorem \ref{mutheorem} below,
that this holds for a pseudo-metric
space as well. Furthermore, Theorem \ref{mutheorem}
contains the corresponding result on the lower
dimension, too, and we now proceed to state the definitions in connection
with this.

\subsection*{The lower dimension of a set}
\begin{definition}
We will say that $(X,d)\in\Lambda_\g$ if there exists
$C(\g)=C(X,d,\g)$ such that, for $x\in X$ and $0<R\le kR\le1$,
$$
N(x,R,k)\ge C(\g)k^\g.\eqno{(\Lambda_\g)}
$$
Then we define the lower dimension $\Lambda(X)$ as:
$$
\Lambda(X)=\sup\{\g,\,(X,d)\in\Lambda_\g\}.
$$
\end{definition}
This dimension was first defined by Larman (\cite{L}) under the name of minimal
dimension. Note that $(X,d)\in\Lambda_0$ is trivial.

\begin{definition}
We will say that a probability measure $\mu$ belongs to $L_\g=L_\g(X,d)$ if
there exists $C(\g)$ such that, for $x\in X$ and $0<R\le kR\le1$,
$$
\mu(B(x,kR))\ge Ck^\g \mu(B(x,R)).\eqno{(L_\g)}
$$
\end{definition}
As before, by taking $k=1/R$ in $(L_\g)$ one gets the weaker condition
$$
\mu(B(x,R))\le C R^\g\qquad x\in X,\quad 0<R\le1.\eqno{(L'_\g)}
$$
Now, observe that the definition
$$
L(X)=\sup\{\g,\,L_\g(X,d)\ne\emptyset\},
$$
will not work. The problem is that $L_\g$
does not imply $\supp(\mu)=X$, so this will say nothing about
$X\setminus \supp(\mu)$. To overcome this problem we make the following
definition.

\begin{definition}
We define the dimension $L(X)$ as:
$$
L(X)=\sup\{\g,\,L_\g(X,d)\cap{\cal U}\ne\emptyset\}.
$$
\end{definition}
Note that $L_0$ poses no restriction on $\mu\in{\cal U}$.

\subsection*{The main theorem}
We now state the main result of this paper. Note that in the special case
when $t=0$, we
can take $t'=t=0$.
\begin{theorem} \label{mutheorem}
Let $(X,d)\in\Upsilon_{s}\cap\Lambda_{t}$, $0\le t\le s<+\infty$, be a
compact complete
pseudo-metric space.
Then for any $s'>s$ and $t'<t$ there exists a probability measure
$\mu\in U_{s'}\cap L_{t'}$.
\end{theorem}
From Theorem \ref{mutheorem} and Propositions \ref{saxo} and
\ref{piano} below we then get

\begin{corolari}
If $\Upsilon(X)<+\infty$, then $\Upsilon(X)=U(X)$ and $\Lambda(X)=L(X)$.
\end{corolari}

\section{Proof of the theorem}
In what follows $X=(X,d)$ denotes an arbitrary compact complete
pseudo-metric space.
To prove Theorem \ref{mutheorem}
we will build a sequence of measures with certain
properties and the measure $\mu$ will be a limit point of this
sequence. We start by proving the trivial
inequalities $\U(X)\le U(X)$ and $\L(X)\ge L(X)$.

\subsection{The trivial inequalities}

\begin{proposition}\label{saxo}
If $\mu\in U_\g(X,d)$, then $X\in \Upsilon_\g$, i.e. $\Upsilon(X)\le U(X)$.
\end{proposition}

\demo{Proof}
Let $\mu\in U_\g$, fix any $x_0\in X$ and let
$x_1,\dots,x_N$ be points in $B(x,kR)$ with $d(x_i,x_j)\ge R$\ for
$i\neq j$. Since $\mu\in U_\g$ and $B(x_0,2C_d kR)\subset B(x_i,4C_d^2 kR)$,
$$
\mu(B(x_0,2C_d kR))\le\mu(B(x_i,4C_d^2 kR))\le C8^\g C_d^{3\g}k^\g
\mu(B(x_i,\frac{R}{2C_d})),
$$
Also, the (open) balls $B(x_i,R/(2C_d))$ are disjoint and lie in
$B(x_0,2C_d kR)$, so
$$
\mu(B(x_0,2C_d kR))\ge\sum_{i=1}^N\mu(B(x_i,\frac{R}{2C_d}))\ge
N\frac{\mu(B(x_0,2C_d kR))}{C8^\g C_d^{3\g}k^\g},
$$
Thus $N\le C(\g,C_d)k^\g$, i.e. $\Upsilon(X)\le U(X)$.\qed

\begin{proposition}\label{piano}
If $\mu\in L_\g\cap{\cal U}$, then $(X,d)\in\Lambda_\g$, i.e.
$\Lambda(X)\le L(X)$.
\end{proposition}

\demo{Proof}
Let $\{x_1,\dots,x_N\}$ be a maximal set of points in
$B(x_0,kR)$ separated by a distance greater than or equal to $R$. Fix any
$\mu\in L_\g\cap{\cal U}$.
Then, since $\mu$ is doubling and $B(x_i,kR)\subset B(x_0,2C_dkR)$ for all $i$,
$$
C_1\mu(B(x_0,kR))\ge\mu(B(x_0,2C_d kR))\ge\mu(B(x_i,kR))\ge Ck^\g
\mu(B(x_i,R)).
$$
Also, since $\{x_1,\dots,x_N\}$ is maximal, $B(x_0,kR)\subset\cup_{i=1}^N
B(x_i,R)$,
$$
\mu(B(x_0,kR))\le\sum_{i=1}^N\mu(B(x_i,R))\le
\frac N{Ck^\g}\mu(B(x_0,kR))
$$
Thus, $N\ge Ck^\g$, i.e. $(X,d)\in \Lambda_\g$.
\qed

\subsection{The main lemma}
Assume that $X=(X,d)\in\Lambda_t\cap\Upsilon_s$, and without loss of
generality suppose
that diam$(X)<1$.
Let $C_d$ be the constant associated to the pseudo-metric $d$,
$C_t$ the constant appearing in $\Lambda_t$ and $C_s$ the one in $\Upsilon_s$.
Given $t'<t$ and $s'>s$, choose $A\ge 16C_d^4$ large enough such that
$A^{s'-s}>C_s$ and
$A^{t-t'}>4^tC_d^{2t}C_t^{-1}$.
For each non-negative integer $j$, let $S_j$ be a maximal set of points in
$X$ separated
by a distance greater than or equal to $A^{-j}$. In particular this means
that $S_0$
consists of just one point.

We define projections $\E=\E_m:S_{m+1}\to S_m$ for $m\ge0$ as follows.
For $g\in S_{m+1}$ choose one of the points $e\in S_m$ for which
$d(g,e)=d(g,S_m)$, and
denote it by $e=\E(g)$.
Then for $e\in S_m$ let
$$
S_{e,m+1}=\{g\in S_{m+1},\,\, e=\E(g)\}.
$$
It is easy to see that $\{S_{e,m+1},\,e\in S_m\}$ is a partition of
 $S_{m+1}$.

The following proposition is a key to the proof of Lemma \ref{lemmavol}
below. The proposition
gives us estimates on the number of points in $S_{e,m+1}$.

\begin{proposition} \label{comptar}
Let $e\in S_m$. then
$$
A^{t'}\le \#(S_{e,m+1})\le A^{s'},
$$
where $\#$ denotes the cardinality of a set.
\end{proposition}

\demo{Proof}
Fix any $e\in S_m$.
Clearly $S_{e,m+1}\subset B(e,A^{-m})$ since $S_m$ is maximal.
Therefore, and since $X\in\Upsilon_s$ and $A^{s'-s}>C_s$,
$$
\#(S_{e,m+1})\le\#(S_{m+1}\cap B(e,A^{-m}))\le N(e,A^{-m-1},A)\le C_sA^s\le
A^{s'},
$$
which proves the right inequality of the proposition.

For the left inequality, we first note that there exists $g\in S_{m+1}$ for
which
$d(g,e)<A^{-m-1}$, and as $A>2C_d$ it is clear that $e=\E(g)$ for such $g$.
Also, for $e'\neq e''$ we have $B(e',A^{-m}/(2C_d))\cap
B(e'',A^{-m}/(2C_d))=\emptyset$.
Thus,
$$
S_{m+1}\cap B(e,A^{-m}/(2C_d))\subset S_{e,m+1}\label{cul}
$$
Next, for $\{g_i\}_{i=1}^n=S_{m+1}\cap B(e,A^{-m}/(2C_d))$ we have
$$
n\ge N(e,A^{-m-1},\frac{A}{2C_d^2}-1)\label{merda}
$$
To check it, suppose the contrary, that is, suppose that
$$
n<N(e,A^{-m-1},\frac{A}{2C_d^2}-1)=n_1.
$$
Then there would exist points $x_1,\dots,x_{n_1}$ in
$B(e,(A/(2C_d^2)-1)A^{-m-1})$
separated from each other
by a distance greater
than or equal to $A^{-m-1}$.

But, for $g\in S_{m+1}\setminus (S_{m+1}\cap B(e,A^{-m}/(2C_d)))$ we have
$$
d(g,x_i)\ge \frac{1}{C_d}d(g,e)-d(e,x_i)\ge \frac{A}{2C_d^2}A^{-m-1}-
\left(\frac{A}{2C_d^2}-1\right)A^{-m-1}=A^{-m-1},
$$
which means that the set
$$
S_{m+1}'=(\{x_i\}_{i=1}^{n_1}\cup S_{m+1})\setminus
(S_{m+1}\cap B(e,\frac{A}{2C_d} A^{-m-1}))
$$
fulfils $\#(S_{m+1}')>\#(S_{m+1})$, a contradiction to the maximality of
$S_{m+1}$.

Thus, from \ref{cul}, \ref{merda},
the choice of $A$ and the fact that $X\in\Lambda_t$,
we conclude
\begin{eqnarray*}
\#(S_{e,m+1})&\ge&\#(S_{m+1}\cap B(e,A^{-m}/(2C_d)))\ge
N(e,A^{-m-1},\frac{A}{2C_d^2}-1)
\ge\\ &\ge& C_t\left(\frac{A}{2C_d^2}-1\right)^t
\ge C_tA^t(4C_d^2)^{-t}\ge A^{t'}
\end{eqnarray*}
\qed

\begin{lemma}\label{lemmavol}
Let $f_0$ be a measure on $S_m$ such that for any $e,e'\in S_m$ we have
$$
f_0(e')\le C_1 f_0(e)
$$
whenever $d(e,e')\le C_2 A^{-m}$, with $C_1=A^{s'-t'}$, and $C_2=8C_d^3$. Then
there is a measure $f_1$ on $S_{m+1}$ with the following properties:
\begin{enumerate}
\item[{\bf (a)}] $f_1(g')\le C_1 f_1(g)$ for any $g,g'\in S_{m+1}$ with
$d(g,g')\le
C_2 A^{-m-1}$.
\item[{\bf (b)}] If $g\in S_{e,m+1}$, then
$A^{-s'}f_0(e)\le f_1(g)\le A^{-t'}f_0(e)$.
\item[{\bf (c)}] $f_0(X)=f_1(X)$.
\item[{\bf (d)}] The construction of the measure $f_1$ from the measure $f_0$
can be regarded as a transfer of mass from the points in  $S_m$
 to those of $S_{m+1}$, with no mass transferred over a distance greater than
$2C_dA^{-m}$. This means that if $g\in S_{m+1}$ recieves mass from
$e\in S_m$, then $d(g,e)\le 2C_d A^{-m}$.
\end{enumerate}
\end{lemma}

\demo{Proof of the lemma}
Let $f_{00}$ be the measure obtained by homogeneously distributing the mass
of each $e\in S_m$ on the points in $S_{e,m+1}$.
By doing so, we obtain a measure satisfying {\bf (b)} (because of
Proposition \ref{comptar}),
{\bf (c)} and {\bf (d)}. If  $f_{00}$ satisfies {\bf (a)}, then let
$f_1=f_{00}$ and we are done.

Assume that $f_{00}$ does not satisfy {\bf (a)}.
Let $\{g_i',g_i''\}_{i=1}^T$ be all the pairs of points in
 $S_{m+1}$ with $d(g_i',g_i'')\le C_2 A^{-m-1}$.
We will construct a finite sequence of measures $\{f_{0j},\,j=1,\dots,T\}$,
such that $f_{0j}$ will satisfy {\bf (a)} for all the pairs
$\{(g_i',g_i'')\}_{i=1}^{j}$, and as we will see $f_1=f_{0T}$ is the
desired measure.

The construction of $f_{0j+1}$ from $f_{0j}$ is as follows:

If $C_1^{-1}f_{0j}(g_{j+1}'')\le f_{0j}(g_{j+1}')\le C_1
f_{0j}(g_{j+1}'')$, then let
$f_{0j+1}=f_{0j}$. Otherwise, only one of these inequalities can fail,
and without loss of generality we may assume that
$f_{0j}(g_{j+1}')> C_1 f_{0j}(g_{j+1}'')$. Then we move mass from $g_{j+1}'$ to
$g_{j+1}''$ by defining $f_{0j+1}$ as
\begin{eqnarray}
f_{0j+1}(g_{j+1}')&=&f_{0j}(g_{j+1}')-
 \frac{f_{0j}(g_{j+1}')-C_1f_{0j}(g_{j+1}'')}{C_1+1};\nonumber\\
f_{0j+1}(g_{j+1}'')&=&f_{0j}(g_{j+1}'')+
 \frac{f_{0j}(g_{j+1}')-C_1 f_{0j}(g_{j+1}'')}{C_1+1};\nonumber\\
f_{0j+1}(g)&=&f_{0j}(g)\qquad\mbox{ if }g\notin
\{g_{j+1}',\,g_{j+1}''\}.\nonumber
\end{eqnarray}
With this definition $f_{0j+1}(g_{j+1}')=C_1f_{0j+1}(g_{j+1}'')$, which
means that
{\bf (a)} is true for $f_{0j+1}$ with respect to $(g_{j+1}',g_{j+1}'')$. In
particular,
note that {\bf (a)} is true for $f_{01}$ with respect to $(g_1',g_1'')$.

We are now going to check condition {\bf (b)} for $f_{0j+1}$.
To do so, suppose that {\bf (b)} holds for $f_{0j}$, i.e. suppose that
$$
A^{-s'}f_0(e)\le f_1(g)\le A^{-t'} f_0(e),\quad g\in S_{e,m+1}.
$$
If $f_{0j+1}=f_{0j}$ or $g\notin\{g_{j+1}',g_{j+1}''\}$, then there is
nothing to check.
Otherwise, as before we can assume that
$f_{0j}(g_{j+1}')> C_1 f_{0j}(g_{j+1}'')$. Let
$e'=\E(g_{j+1}')$ and $e''=\E(g_{j+1}'')$. It is clearly enough to prove that
$f_{0j+1}(g_{j+1}')\ge A^{-s'}f_0(e')$ and
$f_{0j+1}(g_{j+1}'')\le A^{-t'}f_0(e'')$
(because $f_{0j+1}(g_{j+1}')<f_{0j}(g_{j+1}')\le A^{-t'}f_0(e')$ and
$f_{0j+1}(g_{j+1}'')>f_{0j}(g_{j+1}'')\ge A^{-s'}f_0(e'')$).
Now
\begin{eqnarray*}
d(e',e'')&\le& C_d d(e',g_{j+1}')+C_d^2 d(g_{j+1}',g_{j+1}'')+
C_d^2 d(g_{j+1}'',e'')\le\\
&\le&C_d A^{-m}+C_2 C_d^2 A^{-m-1}+C_d^2 A^{-m}\le C_2 A^{-m},
\end{eqnarray*}
so $f_0(e')\le C_1 f_0(e'')$. Therefore
\begin{eqnarray}
f_{0j+1}(g_{j+1}'')=C_1^{-1}f_{0j+1}(g_{j+1}')\le
C_1^{-1}f_{0j}(g_{j+1}')\le C_1^{-1}A^{-t'}f_{0}(e') \le
A^{-t'}f_{0}(e'').\nonumber
\end{eqnarray}
Analogously, $f_0(e'')\ge C_1^{-1}f_0(e')$. Thus,
\begin{eqnarray}
f_{0j+1}(g_{j+1}')=C_1 f_{0j+1}(g_{j+1}'')\ge
C_1 f_{0j}(g_{j+1}'')\ge C_1 A^{-s'}f_{0}(e'')\ge A^{-s'}f_{0}(e').
\nonumber
\end{eqnarray}
Consequently, since {\bf (b)} holds for $f_{00}$ according to
Proposition \ref{comptar} it is
then clear that it holds for $f_1=f_{0T}$ as well.

We are now going to check that when a pair satisfies {\bf (a)} with respect
to $f_{0j}$,
it also does with respect to $f_{0j+1}$. To this end, pick any pair $(g_1,g_2)$,
$d(g_1,g_2)\le C_2 A^{-m-1}$, for which
$$
C_1^{-1}f_{0j}(g_1)\le f_{0j}(g_2)\le C_1 f_{0j}(g_1).
$$
If $(g_1,g_2)$ and $(g_{j+1}',g_{j+1}'')$ have no point in common or
if $f_{0j+1}=f_{0j}$, then we are done. Otherwise,
$f_{0j+1}\ne f_{0j}$ and $f_{0j}(g_{j+1}')> C_1 f_{0j}(g_{j+1}'')$.
Then the two pairs have only one point in common, say $g_1$.
In this case $f_{0j+1}(g_2)=f_{0j}(g_2)$. We have two possible cases to
consider,
either $g_1=g_{j+1}'$ or $g_1=g_{j+1}''$:

If $g_1=g_{j+1}''$, then $f_{0j+1}(g_1)>f_{0j}(g_1)$, so in this case it is
enough to
prove that $f_{0j+1}(g_1)\le C_1 f_{0j+1}(g_2)$. If
$e'=\E(g_{j+1}')$ and $e_2=\E(g_2)$, then
\begin{eqnarray}
d(e',e_2)&\le C_d d(e',g_{j+1}')+C_d^3 d(g_{j+1}',g_1)+
C_d^3 d(g_1,g_2)+C_d^2 d(g_2,e_2)\le\nonumber\\
&\le C_d A^{-m}+2C_d^3C_2 A^{-m-1}+C_d^2 A^{-m}\le C_2 A^{-m},\label{chorra}
\end{eqnarray}
so $f_0(e')\le C_1 f_0(e_2)$.
Also, since we already know that {\bf (b)} is true, we have
$f_0(e_2)\le A^{s'}f_{0j+1}(g_2)$ and $f_{0j+1}(g_{j+1}')\le A^{-t'}f_0(e')$.
Thus,
\begin{eqnarray}
f_{0j+1}(g_1)&=&f_{0j+1}(g_{j+1}'')=C_1^{-1}f_{0j+1}(g_{j+1}')\le
C_1^{-1}A^{-t'}f_{0}(e')
\le\nonumber\\
&\le& A^{-t'}f_{0}(e_2)\le A^{s'-t'}f_{0j}(g_2)=
A^{s'-t'}f_{0j+1}(g_2)=C_1f_{0j+1}(g_2)\nonumber
\end{eqnarray}

Otherwise, if $g_1=g_{j+1}'$, then $f_{0j+1}(g_1)<f_{0j}(g_1)$, and it is
enough to check
that $f_{0j+1}(g_1)\ge C_1^{-1}f_{0j+1}(g_2)$. But, for
$e''=\E(g_{j+1}'')$, then as in \ref{chorra} $d(e'',e_2)\le C_2 A^{-m}$. Also,
$f_{0j+1}(g_1)= C_1 f_{0j+1}(g_{j+1}'')$. Thus, from {\bf (b)} we then get
\begin{eqnarray*}
f_{0j+1}(g_1)&=&C_1 f_{0j+1}(g_{j+1}'') \ge
C_1 A^{-s'}f_{0}(e'')\ge A^{-s'}f_{0}(e_2)\ge\\
&\ge& A^{t'-s'}f_{0j+1}(g_2)=C_1^{-1}f_{0j+1}(g_2).
\end{eqnarray*}

This concludes the proof that {\bf (a)} is true for $f_1$.

Clearly $f_{0j+1}(X)=f_{0j}(X)$, so {\bf (c)} is also true for $f_1$.

It remains to check {\bf (d)}. When passing from
$f_0$ to $f_{00}$ no mass is moved over a distance exceeding $A^{-m}$, because
$S_{e,m+1}\subset B(e,A^{-m})$, and when going from $f_{0j}$ to $f_{0j+1}$
no mass is moved over a distance exceeding $C_2 A^{-m-1}$, and $C_2/A<1$.
It therefore remains to prove that in the construction of $f_1$ from $f_0$
there are no pairs $(g_1,g_2)$ and $(g_2,g_3)$ in $S_{m+1}$ for which mass
is first moved
from $g_1$ to $g_2$ and then at a subsequent step from $g_2$ to $g_3$. To
prove this,
assume the opposite.
Then
$$
f_{00}(g_1)>C_1 f_{00}(g_2)\qquad\mbox{ and }\qquad
f_{00}(g_2)>C_1 f_{00}(g_3).
$$
But, if $e_1=\E(g_1)$ and $e_3=\E(g_3)$, then as in \ref{chorra},
 $d(e_1,e_3)\le C_2 A^{-m}$,
so by the hypothesis $C_1^{-1}f_0(e_1)\le f_0(e_3)\le C_1f_0(e_1)$.
Also,
$$
A^{-s'}f_0(e_i)\le f_{00}(g_i)\le A^{-t'}f_0(e_i),
$$
for $i=1$ and $i=3$. Adding these two inequalities, we would then get
$$
f_0(e_1)\ge A^{t'}f_{00}(g_1)>C_1A^{t'}f_{00}(g_2)>
C_1^2A^{t'}f_{00}(g_3)\ge C_1^2A^{t'-s'}f_0(e_3),
$$
contradicting $f_0(e_1)\le C_1f_0(e_3)$, as $d(e_1,e_3)\le C_2A^{-m}$ and
$C_1=A^{s'-t'}$.
\qed

\subsection{Proof of the theorem}
We will use Lemma \ref{lemmavol}
 to construct a sequence of probability measures and, as we
will see, any limit point of this sequence will satisfy $L_{t'}$ and $U_{s'}$.

We start by defining a probability measure $\mu_0$ on $S_0$ (note that $S_0$
consists of one point only, by the assumption diam$(X)<1$).
Obviously $\mu_0$ satisfies the hypothesis of Lemma \ref{lemmavol}.
For every non-integer $j$ we then use Lemma \ref{lemmavol}
to construct a probability measure
$\mu_{j+1}=f_1$ on $S_{j+1}$ from $\mu_j=f_0$.
In this way we get a sequence $\{\mu_j\}_{j=0}^\infty$ of probability measures.
This sequence lies in the unit ball of the dual of the Banach
space ${\cal C}(X)$, and thus has at least one weak limit point.
Let $\mu$ be any limit point of this sequence. In the proof of the theorem
we will
frequently use the following proposition, based on {\bf (d)} of Lemma
\ref{lemmavol}.

\begin{proposition} \label{cotes}
Let $j\in\nn$, $r\ge 0$ and $x\in X$. Letting $C_4=2C_d^2/(1-C_d/A)$ we
then have
$$
\mu_j(B(x,r))\le \mu(B(x,r+C_4 A^{-j})) %\label{cota_inferior}
$$
and
$$
\mu(B(x,r))\le \mu_j(B(x,r+C_4 A^{-j})). %\label{cota_superior}
$$
\end{proposition}

\demo{Proof}
According to {\bf (d)} of Lemma \ref{lemmavol}
no mass is moved at a distance exceeding $2C_dA^{-j}$ when
constructing
$\mu_{j+1}$ from $\mu_j$. Thus, when passing from $\mu_j$ to $\mu_{j+k}$,
$k\ge1$, no mass is moved at a distance exceeding
$$
2C_d^2 A^{-j}\sum_{n=0}^{k-1}(C_d/A)^{n}<
\frac{2C_d^2}{1-C_d/A}A^{-j}=C_4A^{-j},
$$
which means that there is no mass transfer from $B(x,r)$ into the complement of
$B(x,r+C_4A^{-j})$, and vice versa. Thus,
$$
\mu_j(B(x,R))\le \mu_{j+k}(B(x,r+C_4A^{-j}))
$$
and
$$
\mu_{j+k}(B(x,r))\le\mu_j(B(x,r+C_4A^{-j})).
$$
Now, as $\mu$ is a weak limit point of $\{\mu_{j+k}\}$, the same is true
for $\mu$ as
well.
\qed

We are now going to prove that $\mu\in L_{t'}\cap U_{s'}$. To this end, we
first pick
an $x\in X$, and then some $R$ and $k$ for which $0<R\le kR\le1$. Next we choose
integers $m$ and $M$ such that
\begin{equation}
kR\le A^{-m}< AkR\qquad\mbox{ and }\qquad \frac{R}A\le A^{-M}<R.
\label{kRAkR} \label{indexmM}
\end{equation}
We then denote by $e_{M+1}$ one of the points in $S_{M+1}$ closest to $x$
(there may
be several), and recursively we define $e_{M-j}=\E(e_{M-j+1})\in S_{M-j}$ for
$j=0,\dots,M-m$.

\demo{First claim}
\begin{equation}
\mu_{m+2}(e_{m+2})\le \mu(B(x,kR))\le
C_s3^{s'}(1+C_4)^sC_1\mu_m(e_{m}).\label{tagb}
\end{equation}

\demo{Proof}
To prove this first claim, first note that, by Proposition \ref{cotes},
$$
\mu_{m+2}(e_{m+2})\le \mu(B(e_{m+2},C_4 A^{-m-2})).
$$
On the other hand,
$$
d(x,e_{m+2})\le  C_d A^{-m-2}\sum_{j=0}^{\infty}(C_d/A)^j=
\frac{C_d}{1-C_d/A}A^{-m-2}.
$$
Let $y\in B(e_{m+2},C_4 A^{-m-2})$. Then, by \ref{indexmM},
$$
d(y,x)\le C_dC_4A^{-m-2}+\frac{C_d^2}{1-C_d/A}A^{-m-2}\le A^{-m-1}<kR,
$$
i.e. $B(e_{m+2},C_4A^{-m-2})\subset B(x,kR)$. From Proposition
\ref{cotes} we then get
$$
\mu_{m+2}(e_{m+2})\le\mu(B(e_{m+2},C_4 A^{-m-2}))\le\mu(B(x,kR))
$$
proving the first inequality in \ref{tagb}.
To prove the second inequality, note that \ref{indexmM}
and Proposition \ref{cotes} imply
$$
\mu(B(x,kR))\le\mu_m(B(x,kR+C_4A^{-m}))\le\mu_m(B(x,(1+C_4)A^{-m})).
$$
But, $d(x,e_{m})\le \frac{C_d}{1-C_d/A}A^{-m}$.
Thus, if $e\in S_{m}\cap B(x,(1+C_4)A^{-m})$, then
$$
d(e,e_{m})\le C_d(1+C_4)A^{-m}+\frac{C_d^2}{1-C_d/A}A^{-m}\le
C_2 A^{-m},
$$
so from Lemma \ref{lemmavol} it follows that $\mu_m(e)\le C_1
\mu_m(e_{m})$. Now,
$$
\#\left(S_{m}\cap B(x,(1+C_4)A^{-m})\right)\le C_s(1+C_4)^s,
$$
so from Proposition \ref{cotes} and the fact that $kR\le A^{-m}$, we get
$$
\mu(B(x,kR))\le
\mu_m(B(x,(1+C_4)A^{-m}))\le C_s(1+C_4)^s C_1\mu_m(e_{x,m}),
$$
which concludes the proof of the first claim.

\demo{Second claim}
\begin{equation}
\mu_{M+1}(e_{M+1})\le \mu(B(x,R))\le
C_s(1+C_4)^{s'}A^{2s'}C_1\mu_{M+1}(e_{M+1}).\label{tagc}
\end{equation}

\demo{Proof}
According to Proposition \ref{cotes},
$$
\mu_{M+1}(e_{M+1})\le \mu(B(e_{M+1},C_4A^{-M-1})).
$$
Also, $d(e_{M+1},x)=d(x,S_{M+1})\le A^{-M-1}<R/A$, by the definition of
$e_{M+1}$.
Thus, for $y\in B(e_{x,M+1},C_4A^{-M-1})$,
$$
d(y,x)\le C_dC_4 A^{-M-1}+C_d A^{-M-1}\le A^{-M}<R.
$$
Again by Proposition \ref{cotes},
$$
\mu_{M+1}(e_{M+1})\le\mu(B(e_{M+1},C_4A^{-M-1}))\le\mu(B(x,R)),
$$
proving the left inequality in \ref{tagc}.
To prove the right inequality, note that from Proposition
\ref{cotes} and the fact that
$R\le A^{-M+1}$, by the choice of $M$,
$$
\mu(B(x,R))\le\mu_{M-1}(B(x,R+C_4A^{-M+1}))\le\mu_{M-1}(B(x,(1+C_4)A^{-M+1})).
$$
Also, for $g\in B(x,R+C_4A^{-M+1})\cap S_{M-1}$,
\begin{eqnarray*}
d(g,e_{M-1})\le C_d d(g,x)+C_d^3 d(x,e_{M+1})+
 C_d^3 d(e_{M+1},e_{M})+C_d^2d(e_{M},e_{M-1})\\
\le C_d(1+C_4)A^{-M+1}+C_d^3A^{-M-1}+C_d^3A^{-M}+C_d^2A^{-M+1}\le
C_2A^{-M+1}.
\end{eqnarray*}
Thus, from {\bf (a)} and {\bf (b)} of Lemma \ref{lemmavol} we get (recalling
$e_{M-j}=\E(e_{M-j+1}$),
$$
\mu_{M-1}(g)\le C_1\mu_{M-1}(e_{M-1})\le C_1A^{2s'}\mu_{M+1}(e_{M+1}).
$$
But,
$$
\#\left(B(x,(1+C_4) A^{-M+1})\cap S_{M-1}\right)\le C_s(1+C_4)^s,
$$
so
$$
\mu(B(x,R))\le\mu_{M-1}(B(x,(1+C_4)A^{-M+1}))\le
C_s(1+C_4)^sA^{2s'}C_1\mu_{M+1}(e_{M+1}),
$$
proving the second claim.

To conclude the proof, note that
$\mu(e_{m})\le A^{s'(M+1-m)}\mu_{M+1}(e_{M+1})$ and
$\mu_{m+2}(e_{m+2})\ge A^{t'(M-m-1)}\mu_{M+1}(e_{M+1})$,
by {\bf (b)} in Lemma \ref{lemmavol}.
Also note that $k<A^{M-m}\le A^2 k$, by the choice of $m$ and $M$.
Thus, from the two claims it follows that
$$
\mu(B(x,kR))\le C\mu_m(e_{m})\le
CA^{s'(M-m)}\mu_{M+1}(e_{M+1})\le Ck^{s'}\mu(B(x,R)),
$$
and similarly,
$$
\mu(B(x,kR))\ge\mu_{m+2}(e_{m+2})\ge CA^{t'(M-m)}\mu_{M+1}(e_{M+1})\ge
Ck^{t'}\mu(B(x,R)),
$$
i.e. $\mu\in\Lambda_{t'}\cap\Upsilon_{s'}$.

Note that the final constants $C$ depend only on the given constants $C_d$,
$C_s$, $C_t$
and the choice of $A$, $s'$ and $t'$.
Also note that the last inequality depends on the fact that
$\Upsilon(X)<+\infty$.
\qed

\section{Examples}

Denoting by $\dim(E)$ the Hausdorff dimension of $E$, Larman (\cite{L})
proved that
$\Lambda(E)\le \dim(E)\le \Upsilon(E)$ for any metric space $E=(E,d)$.

We consider some examples of Cantor type sets in the Euclidean metric.

When $E$ is the usual $s$-dimensional Cantor set the dimensions coincide,
$\Lambda(E)=\dim(E)=\Upsilon(E)$ (since $H^s(B(x,r))\approx r^s$ for the
$s$-dimensional Hausdorff measure $H^s$).

Now, denote by $C_t$ the $t$-dimensional Cantor set contained in $[0,1]$
and by $C_s$
the $s$-dimensional Cantor set in $[2,3]$, and suppose that $0<t<s<1$.
Then $E=C_t\cup C_s$ has lower dimension $t=\Lambda(E)=\dim(C_t)<\dim(E)$
and upper
dimension $s=\Upsilon(E)=\dim(C_s)=\dim(E)$ (which follows from the fact that
$\mu(E)=H^t(C_t)+H^s(C_s)\in U_s\cap L_t$,
i.e  $t\le\Lambda(E)\le\Upsilon(E)\le s$, and also,
$s=\dim(E)\le\Upsilon(E)$ and
$t=\Lambda(C_t)\ge \Lambda(E)$).

Next we consider the case when two Cantor sets intersect at the endpoints.
Let $C_1$ be the $\log2/\log3$-dimensional Cantor set in
$[0,1]$ and $C_2$ the $\log2/\log9$-dimensional Cantor set in $[1,2]$
(i.e. obtained by indefinitely deleting $7/9$ from the middle of each
sub-interval
starting with $[1,2]$). It is clear $C_1$ and $C_2$ are closed,
and that $C_1\cap C_2=\{1\}$.
Put $F=C_1\cup C_2$. Then $F$ is closed, and it is easy to see that $F$ has
lower
dimension $\Lambda(F)=\log2/\log9$ and upper dimension
$\Upsilon(F)=\log2/\log3$.
Thus, it follows from Theorem 2
that there exists $\mu\in L_t\cup U_s$ on $F$ for every $t<\log2/\log9$
and $s>\log2/\log3$.
Furthermore, there is a $\mu_n\in L_{\log2/\log9}\cup U_{s_n}$ on $F$ with
$\log2/\log3<s_n<\log2/\log3+1/n$, for any $n\in\Bbb N$.
To see this, choose the constant $A$ in Theorem 2 on the form $A=9^n$,
$n\in\Bbb N$.
Then the corresponding maximal $9^{-mn}$-sets $S_m$ turns out to be
uniquely determined and $2^n\le\#S_{e,m+1}\le4^n+2^n-1$ for all $e$ and $m$
(namely $\#S_{e,m+1}=4^n$ for $e\in[0,1)$, $\#S_{e,m+1}=4^n+2^n-1$ for
$e=1$, and
$\#S_{e,m+1}=2^n$ for $e\in(1,2]$).
This bound on $\#S_{e,m+1}$ imply that there is a  measure
$\mu_n\in L_{\log2/\log9}\cup U_{s_n}$ on $F$,
with $s_n=\log(4^n+2^n-1)/\log9^n<\log2/\log3+1/n$.

Similarly, there is a measure $\mu\in L_t\cap U_{\log2/\log3}$ on $F$ for some
$t<\frac{\log2}{\log9}$.
Let $\nu_1$ be any $\log2/\log3$-measure on $C_1$ and
$\nu_2$ any $\log2/\log9$-measure on $C_2$. Then one can easily see that
the measure $\nu=\nu_1+\nu_2$ is not even doubling on $F$
(because of the density jump at $\{1\}$).
But, it is possible to construct a measure $\mu_2$ on $C_2$ such that
$\mu=\nu_1+\mu_2$ becomes a measure on $F$ belonging to $L_t\cap
U_{\log2/\log3}$ for some
$t<\log2/\log9$ (for details, see \cite{Per}, Part C, Ex. 2).

\smallskip

We do not know whether there exists a
$\mu\in L_{\log2/\log9}\cap U_{\log2/\log3}$ on $F$.
\medskip

It is important to notice that the measure we construct in the proof of the
theorem \ref{mutheorem} is not unique. This was already noticed by Volberg
and Konyagin (\cite{V-K}) and in fact, as proved in \cite{Kaufman-Wu},
in any closed perfect set there are at least two mutually singular doubling
measures, and as they are doubling both of them must lie in some $U_\g$.

However, in spite of this, one can show that Besov spaces
defined with respect to different such measures are equivalent (see
\cite{Per}, Part B).

\end{document}